\documentclass[11pt]{article}
\usepackage{latexsym,amsmath,amsfonts,amssymb}
\def\qed{\hfill$\Box$}

\newcommand{\leqa}{\mbox{$ \;\stackrel{(a)}{\leq}\; $}}

\newcommand{\eqa}{\mbox{$ \;\stackrel{(a)}{=}\; $}}
\newcommand{\eqb}{\mbox{$ \;\stackrel{(b)}{=}\; $}}

\newcommand{\RL}{{\mathbb R}}

\newcommand{\IN}{{\mathbb Z}}

\newcommand{\PR}{\mbox{\rm Pr}}

\newcommand{\la}{\lambda}

\newcommand{\Pol}{\mbox{\rm Po}(\la)} 
\newcommand{\Pot}{\mbox{\rm Po}(t)}

\newcommand{\stv}{{\mbox{\scriptsize TV}}}

\newcommand{\iid}{\mbox{i.i.d.}\!}

\def\be{\begin{eqnarray}}
\def\ee{\end{eqnarray}}
\def\ben{\begin{eqnarray*}}
\def\een{\end{eqnarray*}}

\newcommand{\ro}{\rho}

\newlength{\noteWidth}
\long\def\notes#1{\ifinner
             {\tiny #1}
             \else
             \marginpar{\parbox[t]{\noteWidth}{\raggedright\tiny #1}}
             \fi}

\title{Entropy and the Law of Small Numbers}

\author{
I. Kontoyiannis\thanks{ 
	Division of Applied Mathematics 
	and Dept of Computer Science, 
	Brown Univ, 182 George St., 
	Providence, RI 02912, USA. 
	Email: {\tt yiannis@dam.brown.edu}
	Web: {\tt www.dam.brown.edu/people/yiannis/}.
	Supported in part by NSF grants \#0073378-CCR
	and DMS-9615444, and by USDA-IFAFS grant 
	\#00-52100-9615.
		       }
\and
P. Harremo\"{e}s\thanks{
	Department of Mathematics,
	University of Copenhagen,
	Universitetsparken 5,
 	DK-2100 K{\o}benhavn $\emptyset$,
	Denmark.
	Email: {\tt moes@math.ku.dk}.
	Supported in part by a grant from the Cowi Foundation,
	and by a post-doctoral fellowship
	from the Villum Kann Rasmussen Foundation.
		       }
\and
O. Johnson\thanks{Statistical Laboratory, Centre for 
	Mathematical Sciences,
	Wilberforce Road, Cambridge CB3 0WB, U.K.
	 Email: {\tt otj1000@cam.ac.uk}.
		 }
}
 
\date{\today}
\begin{document}
\bibliographystyle{plain}
\maketitle
\begin{abstract}
Two new information-theoretic methods are introduced
for establishing Poisson approximation inequalities. 
First, using only elementary information-theoretic techniques
it is shown that,
when $S_n=\sum_{i=1}^nX_i$ 
is the sum of the (possibly dependent) binary
random variables $X_1,X_2,\ldots,X_n$,
with $E(X_i)=p_i$ and $E(S_n)=\la$, then
\ben
D(P_{S_n}\|\Pol)\leq
\sum_{i=1}^n p_i^2 \;+\;
\Big[\sum_{i=1}^nH(X_i) \,-\,H(X_1,X_2,\ldots, X_n)\Big],
\een
where $D(P_{S_n}\|\mbox{Po}(\la))$ is the
relative entropy between the distribution 
of $S_n$ and the Poisson($\la$)
distribution. The first term in this bound
measures the individual smallness of the $X_i$
and the second term measures their dependence. 
A general method is outlined for 
obtaining corresponding bounds 
when approximating the distribution 
of a sum of general discrete random 
variables by an infinitely divisible 
distribution.

Second, in the particular case when the $X_i$ 
are independent, the following
sharper bound is established,
\ben
D(P_{S_n}\|\Pol)\leq
\frac{1}{\lambda} \sum_{i=1}^n \frac{p_i^3}{1-p_i},
\een
and it is also generalized to the case
when the $X_i$ are general integer-valued 
random variables. Its proof is based on the 
derivation of a subadditivity property
for a new discrete version of the 
Fisher information, and uses a recent 
logarithmic Sobolev inequality for the 
Poisson distribution.

\end{abstract}
 
{\bf Keywords: } Poisson approximation,
law of small numbers, convergence in relative entropy,
Fisher information, total variation, logarithmic Sobolev
inequality, subadditivity

\newpage
\section{Introduction}

Let $X_1,X_2,\ldots,X_n$ be 
binary random variables. A classical 
result in probability states that, 
if the $X_i$ are independent and
identically distributed ($\iid$)
with common parameter $p_i=E(X_i)=\la/n$, 
then, when $n$ is large, the distribution of their sum
$$S_n=X_1+X_2+\cdots+X_n$$
is close to $\Pol$,
the Poisson distribution
with parameter $\lambda$.
More generally, analogous results
apply when the $X_i$ 
are possibly dependent and not
necessarily identically distributed.  
The distribution of $S_n$ is close
to $\Pol$ as long as:
\begin{itemize}
\item[$(a)$]
The sum $\sum p_i$ of the parameters $p_i$ 
of the $X_i$ is close to $\lambda$.
\item[$(b)$]
None of the $X_i$ dominate the sum, i.e.,
all the $p_i$ are small.
\item[$(c)$]
The variables $X_i$ are not strongly dependent.
\end{itemize}
Such results are often referred to as
``laws of small numbers'' or
``Poisson approximation results.'' 
See 
\cite{arratia-g-g:89}\cite[Section~2.6]{durrett:book}\cite{barbour-et-al:book} 
for details.

Our purpose here is to illustrate 
how techniques based on information-theoretic
ideas can be used to establish
general Poisson approximation inequalities.
In Section~\ref{s:HK} we prove:

\medskip

\noindent
{\em Proposition 1. Poisson Approximation in Relative Entropy}:
If $S_n=\sum_{i=1}^nX_i$ is the sum of~$n$ (possibly
dependent) binary random variables $X_1,X_2,\ldots,X_n$
with parameters $p_i=E(X_i)$ and
with $E(S_n)=\sum_{i=1}^n p_i=\la,$
then the 
distribution $P_{S_n}$ of $S_n$ 
satisfies
\be
D(P_{S_n}\|\Pol)\leq
\sum_{i=1}^n p_i^2 \;+\;
\Big[\sum_{i=1}^nH(X_i) \,-\,H(X_1,X_2,\ldots,X_n)\Big].
\label{eq:intro}
\ee

For two probability distributions 
$P$ and $Q$ on a discrete set $S$,
the relative entropy 
between $P$ and $Q$ is defined as
$D(P\|Q)=\sum_{x\in S} P(x)\log\frac{P(x)}{Q(x)},$
and the entropy of a discrete random
variable (or random vector) $X$ with 
distribution $P$ on $S$ is 
$H(X)=H(P)=-\sum_{x\in S}P(x)\log P(x),$
where $\log$ denotes the natural logarithm.

Whenever $(a)$, $(b)$ and~$(c)$ 
hold we expect the two terms in the 
right-hand side of (\ref{eq:intro}) 
to be small, and hence the distribution 
of $S_n$ to be close to $\Pol$ in the relative
entropy sense. Although
$D(P\|Q)$ is not a proper metric,
it is a natural
measure of ``dissimilarity'' in the context
of statistics 
\cite{kullback-book}\cite[Ch.~12]{cover:book},
and it can be used to define a topology 
on probability measures
\cite{harremoes:topology}.
Also, bounds in relative entropy can
be translated into bounds in total variation
via Pinsker's inequality \cite{cover:book}
\be
\frac{1}{2}\left\|P-Q\right\|_\stv^2\leq 
D(P\|Q).
\label{eq:pinsker}
\ee
For example, if the $X_i$ are independent
(\ref{eq:intro}) reduces to
\be
D(P_{S_n}\|\Pol)\leq\sum_{i=1}^n p_i^2.
\label{eq:simple}
\ee
Although this is reminiscent
of the simple total-variation
bound due to Le Cam \cite{lecam:60},
$$
\left\|P_{S_n}-\Pol\right\|_\stv\leq
\sum_{i=1}^n p_i^2$$
(which, incidentally, only holds when the $X_i$
are independent),
applying Pinsker's inequality (\ref{eq:pinsker})
to (\ref{eq:simple}) leads to the 
suboptimal bound
\be
\left\|P_{S_n}-\Pol\right\|_\stv\leq
\Big[2\sum_{i=1}^n p_i^2\Big]^{1/2}.
\label{eq:rough}
\ee

The proof of Proposition~1 uses only 
elementary information-theoretic facts
that are established using little
more than Jensen's inequality. 
To get sharper bounds for the 
case of independent random variables $X_i$,
in Section~\ref{s:tv} we 
employ a new discrete
version of the Fisher information which we
call {\em scaled Fisher information},
and we prove:

\smallskip

\noindent
{\em Theorem 1. Poisson Approximation for Independent Variables}:
If $S_n=\sum_{i=1}^nX_i$ is the sum of $n$ 
independent binary random variables $X_1,X_2,\ldots,X_n$,
with $E(S_n)=\sum_{i=1}^n p_i=\la,$ 
\be
D(P_{S_n}\|\Pol)\leq
\frac{1}{\lambda} \sum_{i=1}^n \frac{p_i^3}{1-p_i}.
\label{eq:intro2}
\ee

The proof of Theorem~1 combines 
a natural discrete analog of Stam's subbativity of 
the Fisher information \cite{stam:59}\cite{blachman:65}, 
and a recent logarithmic Sobolev inequality of Bobkov and Ledoux 
\cite{bobkov-ledoux:98}.
As we discuss extensively in Section~\ref{s:tv}, 
Theorem~1 is a significant improvement 
over Proposition~1, and in certain cases 
it leads to total variation bounds
that are asymptotically optimal up to 
multiplicative constants in the
convergence rate. Moreover, (\ref{eq:intro2})
is a nontrivial improvement over existing results, 
as it gives a bound for the relative entropy 
and not just the total variation distance.

For an information-theoretic interpretation,
consider a triangular array of 
binary random variables 
$\{(X^{(n)}_1,X^{(n)}_2,\ldots,X^{(n)}_n),\; n\geq 1\}$,
such that the right-hand side of (\ref{eq:intro}) 
goes to zero as $n\to\infty$ (as, for example,
when the $X^{(n)}_i$ 
are $\iid\,$ Bernoulli($\la/n$)). Then 
the distribution of $S_n$ 
converges to $\Pol$, i.e.,
$P_{S_n}$ comes closer 
and closer to the ``most random'' 
distribution among all those that can 
be obtained by summing a finite number 
of Bernoulli random variables:
Let ${\cal P}(\la)$ denote the set of 
all distributions of sums $S_n$ of $n$ 
{\em independent} binary random variables 
with $E(S_n)=\la$, for any finite $n$.
Then \cite{harremoes:01},
\ben
H(\Pol)
\;=\;
\sup \{H(P)\;:\;P\in{\cal P}(\la)\}.
\een
So, roughly and somewhat incorrectly speaking,
the entropy of $S_n$ ``increases'' to the maximum 
entropy $H(\Pol)$
as $n$ grows.
This invites a tempting analogy 
with the second law of thermodynamics, 
stating that the uncertainty of a physical 
system increases with time, until the system 
reaches equilibrium in its maximum entropy state.

Corresponding information-theoretic interpretations 
and proofs have been given for
numerous classical results of proability theory,
including the central limit theorem 
\cite{linnik:59}\cite{brown:82}\cite{barron:clt}\cite{johnson:00},
the convergence of Markov chains 
\cite{renyi:61}\cite{kendall:63}\cite{barron:isit00}, 
many large deviations results
\cite{csiszar:84}\cite{dembo-zeitouni:book}%
\cite{csiszar-cover-choi:87},
the martingale convergence theorem
\cite{barron:isit91}\cite{barron:isit00},
and the Hewitt-Savage 0-1 law \cite{oconnell:TR}.
See also the powerful comments in
\cite[pp.~211,215]{genedenko-korolev:book}.
Finally, we mention that Johnstone and MacGibbon 
considered the problem of Poisson convergence 
from the information theory angle
in \cite{johnstone-macginbbon:87}.
Their approach is different from ours,
and parallels that in \cite{brown:82}\cite{barron:clt} 
for the central limit theorem.


\section{General Bounds in Relative Entropy}
\label{s:HK}


Before giving the proof of Proposition~1 we introduce
some some notation and briefly recall two 
elementary, well-known facts.
The first one formalizes the intuitive idea that 
we cannot do better in a hypothesis test by simply 
pre-processing the data.
Suppose $X$ and $Y$ are random variables
with distributions $P$ and $Q$, 
respectively, let $f$ be an
arbitrary function, and write $P', Q'$ 
for the distribution
of $f(X)$ and $f(Y)$, respectively.
The following ``data processing''
inequality is an easy consequence 
of Jensen's inequality
\cite[Lemma~1.3.11]{csiszar:book},
$$D(P'\|Q')\leq D(P\|Q).$$
Next, given $X$ and $Y$ with
joint distribution $P_{X,Y}$ and 
marginals $P_X$ and $P_Y$, let
$I(X;Y)= H(X)-H(X|Y)$ denote
their mutual information.
The ``chain rule''
is the simple expansion,
$$D(P_{X,Y}\|Q_X\times Q_Y)=D(P_X\|Q_X)+D(P_Y\|Q_Y)+I(X;Y),$$
for any two probability distributions
$Q_X$ and $Q_Y$.

\smallskip

%
%

\noindent
{\em Proof of Proposition 1. }
If we define $S_n'=\sum_{i=1}^n Z_i$,
where $Z_i$ are independent Poisson($p_i$) random
variables, then the distribution 
$P_{S'_n}$ of $S_n'$ is $\Pol$ and
\be
D(P_{S_n}\|\Pol)
&=&
D(P_{S_n}\|P_{S'_n})
	\nonumber\\
&\leqa&
D(P_{X_1,\ldots,X_n}\|P_{Z_1,\ldots,Z_n})
	\nonumber\\
&\eqb&
\sum_{i=1}^n D(P_{X_i}\|\mbox{Po}(p_i))
+\sum_{i=1}^{n-1} I(X_i;(X_{i+1},\ldots,X_n)),
\label{eq:ex1}
\ee
where $(a)$ follows from the data processing 
inequality, and $(b)$ follows by applying 
the chain rule $(n-1)$ times. 
Using simple calculus we obtain the bound 
$$
D(\mbox{Bern}(p)\|\mbox{Po}(p))
=(1-p)\log\frac{(1-p)}{e^{-p}}
	+ p\log\frac{p}{pe^{-p}}\leq p^2,$$
which, applied to each term in the first sum in (\ref{eq:ex1}),
gives,
\be
D(P_{S_n}\|\Pol)
&\leq&
	\sum_{i=1}^n p_i^2 + \sum_{i=1}^{n-1} I(X_i;(X_{i+1},\ldots,X_n))
	\label{eq:ex1b}\\
&=&
	\sum_{i=1}^n p_i^2 + 
	\Big[\sum_{i=1}^n H(X_i) \,-\,H(X_1,X_2,\ldots,X_n)\Big],
	\nonumber
\ee
where in the last step we expanded 
the definition of the mutual informations.
\qed


\newpage


The first term in the above bound makes precise 
what we mean by the requirement that ``all the $p_i$ be small''
whereas the second term quantifies their degree of dependence.
It is worth noting that this difference between the sum of the
entropies of the $X_i$ and their joint entropy can also be written
as the relative entropy $D(P_{X_1^n}\|P_{X_1}\times\cdots\times P_{X_n})$
between their joint distribution and the product
of their marginals. This expression also admits a natural
interpretation as a measure of how far the $X_i$ are from
being independent.

As indicated in the introduction, although the result of 
Proposition~1 is generally good enough to prove convergence
to the Poisson distribution, for finite $n$ it often
gives a suboptimal convergence rate. This is also 
illustrated in the following two examples.

\medskip

\noindent
{\em A Markov Chain.}
Let
$\{(X^{(n)}_1,X^{(n)}_2,\ldots,X^{(n)}_{n}),\;\; n\geq 1\}$
be a triangular array of binary random variables
such that each row
$(X^{(n)}_1,\ldots,X^{(n)}_{n})$
is a Markov chain with
transition matrix
$$
\left(
  \begin{array}{cc}
  \frac{n}{n+1} & \;\;\frac{1}{n+1}\\
  \;		  & \;		     \\
  \frac{n-1}{n+1} & \;\;\frac{2}{n+1}\\
  \end{array}
\right)
$$
and with each $X^{(n)}_i$ having 
(the stationary) Bernoulli($\frac{1}{n}$)
distribution. The convergence of the 
distribution of 
$S_n=\sum_{i=1}^nX^{(n)}_i$
to Po(1) 
is a well-studied problem;
see 
\cite{cekanavicius:02}
and the references therein.
Applying Proposition~1
(or, equivalently, inequality (\ref{eq:ex1b}))
in this case translates to
$$
D(P_{S_n}\|\mbox{Po}(1))
\;\leq\;
	\sum_{i=1}^n\frac{1}{n^2}
	+\sum_{i=1}^{n-1}I(X^{(n)}_i;X^{(n)}_{i+1})
\;=\;
	\frac{1}{n} + (n-1)I(X^{(n)}_1;X^{(n)}_2),
$$
since 
$I(X^{(n)}_i;(X^{(n)}_{i+1},\ldots,X^{(n)}_n))
=I(X^{(n)}_i;X^{(n)}_{i+1})$ by
the Markov property,
and stationarity implies that 
$I(X^{(n)}_i;X^{(n)}_{i+1})=I(X^{(n)}_1;X^{(n)}_2)$.
A straightforward calculation yields that
$$(n-1)I(X^{(n)}_1;X^{(n)}_2)
=
(n-1)
	\left[h({\textstyle\frac{1}{n}})
	-h({\textstyle\frac{1}{n+1}})\right]
+\frac{n-1}{n}h({\textstyle\frac{1}{n+1}})
-\frac{n-1}{n}h({\textstyle\frac{2}{n+1}}),
$$
where $h(p)$ denotes the binary entropy
function $h(p)=-p\log p -(1-p)\log(1-p)$,
and simple calculus shows that all three
terms above converge to zero as $n\to\infty$.
In fact, this expression can be bounded above
by
$$h({\textstyle\frac{1}{n+1}})+\frac{\log n}{n}
\leq 3\frac{\log n}{n},$$
where the last inequality holds for all
$n\geq 3$,
so putting it all together,
$$D(P_{S_n}\|\mbox{Po}(1))
\;\leq\;3\frac{\log n}{n} + \frac{1}{n}.$$
[A corresponding bound can similarly be derived 
if instead of stationarity we assume
that $X^{(n)}_1$ has 
$p^{(n)}_1=E(X^{(n)}_1)<1/n$.]
As mentioned above, although this bound is sufficient
to prove that $P_{S_n}$ converges to the Poisson 
distribution, it leads to a convergence rate in
total variation of order $\sqrt{(\log n)/n}$,
compared to the $O(1/n)$ bound derived
in \cite{barbour-et-al:book}\cite{serfozo:86}\cite{serfozo:88}.

\medskip

\noindent
{\em A Compound Poisson Approximation Example. }
Let $X_1,\ldots,X_n$ be independent
Bernoulli random variables with parameters
$p_i=E(X_i)$, write $\la=\sum_{i=1}^n p_i$,
and let $\alpha_1,\alpha_2,\ldots,
\alpha_n$ be $\iid$, independent of the $X_i$,
with distribution
$$
\alpha_i=\left\{
		\begin{array}{cc}
  			1 & \mbox{with prob}\;1/2\\
  			2 & \mbox{with prob}\;1/2.
  \end{array}
\right.
$$
We will show that the distribution of the
sum 
$$S_n=\sum_{i=1}^n\alpha_i X_i$$
is close to
the compound Poisson distribution with 
parameters
($\la/2,\la/2$), which we denote by
Po($\la/2,\la/2$).
Recall that if $Z_1$ and $Z_2$ are $\iid$ Poisson($\la/2$)
random variables, then $Z=(Z_1+2Z_2)$
has $\mbox{Po}(\la/2,\la/2)$ distribution.
Alternatively, we can write
$Z=\sum_{i=1}^nY_i$ where the $Y_i$ are
independent Po($p_i/2,p_i/2$) random 
variables. Arguing as before,
the data processing inequality and
the chain rule imply that
$$D(P_{S_n}\|\mbox{Po}(\la/2,\la/2))
\leq
	D(P_{\alpha_1X_1,\ldots,\alpha_nX_n}
	\|P_{Y_1,\ldots,Y_n})
=
\sum_{i=1}^nD(P_{\alpha_iX_i}\|P_{Y_i}),$$
and it is straightforward to calculate
$$
D(P_{\alpha_iX_i}\|P_{Y_i})
\leq
	p_i^2 + (1-p_i)[p_i+\log(1-p_i)]
	-\frac{p_i}{2}\log(1+p_i/4)
\leq
	p_i^2,
$$
so that
$$D(P_{S_n}\|\mbox{Po}(\la/2,\la/2))\leq \sum_{i=1}^np_i^2.$$

\noindent
{\em A general method. }
Finally, we outline a simple general strategy
for approximating
the distribution $P_{S_n}$ of the sum of
$n$ nonnegative-integer-valued random variables
$X_1,X_2,\ldots,X_n$ by
the distribution of 
some infinitely divisible
discrete
random variable $Z$ with
$E(S_n)=E(Z)$.

First, use the infinitely divisibility of $P_Z$ 
to represent $Z$ as $Z=\sum_{i=1}^nY_i$
where the $Y_i$ are independent and have the 
same distribution as $Z$ but with different 
parameters. Then apply the data processing 
inequality and the chain rule as before 
to obtain
$$D(P_{S_n}\|P_Z)\leq
        \sum_{i=1}^nD(P_{X_i}\|P_{Y_i})
	+ \Big[\sum_{i=1}^nH(X_i) \,-\,H(X_1,\ldots X_n)\Big],$$
and finally, estimate the last two terms in above inequality.
The first term should be small if the
$X_i$ are individually small and well-approximated
by the corresponding $Y_i$, and the second term
should be small if the $X_i$ are sufficiently 
weakly dependent.


\section{Tighter Bounds for Independent Random Variables}
\label{s:tv}


Next we take a different point of view that
yields tighter bounds than Proposition~1. 
Recall that in
\cite{johnstone-macginbbon:87}\cite{papathanasiou:93}\cite{kagan:01},
the Fisher information
of a random variable $X$ with distribution
$P$ on $\IN_+=\{0,1,2,\ldots\}$, 
is defined in a way analogous to that for 
continuous random variables, via
\ben
J(X) = E\Big[\Big(\frac{P(X-1)-P(X)}{P(X)}\Big)^2\Big],
\een
with the convention that $P(-1)=0$.
However, as Kagan \cite{kagan:01} acknowledges,
this definition is really only useful if $X$ is
supported on the entire $\IN_+$:
If $X$ has bounded support then for some $n$, $P(n) > 0$ 
but $P(n+1) = 0$, which implies that $J(X)=\infty.$


Partly in order to avoid this difficulty, we proceed along
a different route. Recalling that the Poisson distribution
is characterized by the recurrence $\la P(x)=(x+1) P(x+1)$ for all $x$,
we let the {\em scaled score function} of a random
variable $X$ with mean $\la$ and distribution $P$ on $\IN_+$ be
$$\ro_X (x) = \frac{ (x+1)P(x+1)}{\la P(x)} - 1,
\;\;\;\;x\in\IN_+,$$
and  
we define the {\em scaled Fisher information of $X$} as
\ben
K(X)= \la E[\rho_X(X)^2].
\een
From this we easily see that 
$$K(X)\geq 0$$
with equality iff $\ro_X (X)=0$ with probability 1,
i.e., iff $X$ is has a Poisson$(\la)$ distribution. Moreover,
as we show next, the smaller
the value of $K(X)$, the closer $P$ is to the Poisson($\la$)
distribution. The proof of Proposition~2, given in Section 3.2,
is an easy consequence of a recent logarithmic Sobolev 
inequality of Bobkov and Ledoux \cite{bobkov-ledoux:98}.

\medskip

\noindent
{\em Proposition 2. Relative Entropy and $K(X)$}:
If $X$ is a random variable with distribution $P$ on $\IN_+$
and with $E(X)=\la$, then
\be
D(P\|\Pol)\leq K(X),
\label{eq:ls}
\ee
as long as either $P$ has full support (i.e., $P(k)>0$ for all $k$), 
or finite support (i.e., there exists $N\in\IN_+$ such that $P(k)=0$
for all $k>N$).

\medskip

Note that from (\ref{eq:ls}) and Pinsker's inequality 
(\ref{eq:pinsker}) we have that
\be
\|P-\Pol\|_\stv\leq\sqrt{2 K(X)}.
\label{eq:tvbound}
\ee
We also give a direct proof of (\ref{eq:tvbound})
in Section~3.2, based on a simple
Poincar\'{e} inequality for the Poisson measure.

\subsection{Results}

The main step in the proof of Theorem~1
will be to establish a form of
subadditivity for the scaled Fisher information.
It is worth noting that
in the Gaussian case the Fisher information 
is also subadditive \cite{stam:59}\cite{blachman:65},
but, in contrast to the present setting,
subadditivity alone does not suffice
to prove the central limit theorem \cite{barron:clt}.
Proposition~3 is proved in Section~3.2.

\medskip

\noindent
{\em Proposition 3. Subadditivity of Scaled Fisher Information}:
If $S_n=\sum_{i=1}^nX_i$ is the sum of $n$
independent integer-valued random variables 
$X_1,X_2,\ldots,X_n$, with means $E(X_i)=p_i$
and $E(S_n)=\sum_{i=1}^n p_i=\la,$ then
\ben
K(S_n) \leq \sum_{i=1}^n \frac{p_i}{\la} K(X_i).
\een

\noindent
{\em Proof of Theorem 1. }
If the $X_i$ are independent Bernoulli($p_i$) random 
variables with $\sum_{i=1}^np_i=\la$, then
$K(X_i) = p_i^2/(1-p_i)$ and Proposition~3 gives
$$ K(S_n) \leq \frac{1}{\lambda} 
\sum_{i=1}^n \frac{p_i^3}{1-p_i}.$$
Combining this with $X=S_n$ in Proposition~2 yields
inequality (\ref{eq:intro2}).\qed

\medskip

\noindent
{\em Example 1. }
If the $X_i$ are $\iid$ Bernoulli($\la/n$)
random variables, from Theorem~1 combined
with Pinsker's inequality (\ref{eq:pinsker})
we obtain that for any $\epsilon>0$,
$$\left\|P_{S_n}-\mbox{Po}(\lambda)\right\|_\stv \leq
(2+\epsilon)\frac{\la}{n},\;\;\;\;\mbox{ for $n\geq\la/\epsilon$.}$$
This is a definite improvement over the earlier
$2\la/\sqrt{n}$ bound from (\ref{eq:rough}),
and, except for the constant factor,
it is asymptotically of the right order;
see \cite{barbour-et-al:book}\cite{deheuvels-pfeifer:86}
for details. 


\medskip

\noindent
{\em Example 2. }
If the $X_i$ are $\iid$ Bernoulli($\mu/\sqrt{n}$)
random variables, Theorem~1 together with
Pinsker's inequality (\ref{eq:pinsker}) yield,
$$\left\|P_{S_n}-\mbox{Po}(\mu\sqrt{n})\right\|_\stv \leq
\frac{\mu}{\sqrt{n}}\sqrt{
\frac{2}{1-\mu/\sqrt{n}}
}\approx \frac{\mu}{\sqrt{n}}\sqrt{2},
$$
which is of the same order as the optimal asymptotic rate,
as $n\to\infty$,
$$\left\|P_{S_n}-\mbox{Po}(\mu\sqrt{n})\right\|_\stv 
\sim\frac{\mu}{\sqrt{n}}\sqrt{1/(2\pi e)}$$
derived in \cite{deheuvels-pfeifer:86}.

\medskip

\noindent
{\em Example 3. }
If the $X_i$ are Geometric random variables
with respective distributions $P_i(x)=(1-q_i)^xq_i,$
$x\geq 0$,
then $K(X_i)=(1-q_i)^2/q_i$.
Letting $S_n=\sum_{i=1}^nX_i$
and assuming that $E(S_n)=\sum_{i=1}^n\frac{1-q_i}{q_i}=\la$,
combining Proposition~3 and the bound (\ref{eq:tvbound})
yields
$$\left\|P_{S_n}-\Pol\right\|_\stv\leq\sqrt{\frac{2}{\la}\sum_{i=1}^n
\frac{(1-q_i)3}{q_i^2}}.$$
In particular, taking all the $q_i=n/(n+\la)$ gives
the elegant estimate
$$\left\|P_{S_n}-\Pol\right\|_\stv\leq\frac{\sqrt{2}\la}{\sqrt{n(n+\la)}}
\leq\sqrt{2}\frac{\la}{n}.$$

\medskip

To see how tight the result of Proposition~3 is
in general, note that 
the following lower bound of Cram\'{e}r-Rao type
holds: Since for all $a$ and any random variable $S$ with
mean $\la$ and variance $\sigma^2$,
\be
0 \leq \la E( \ro_S(S) - a(S-\la))^2 = K(S) + \la \left( a^2 \sigma^2 - 2a 
\left( \frac{ \sigma^2 - \la}{\la} \right) \right),
\label{eq:CR}
\ee
choosing $a = (\sigma^2 - \la)/(\sigma^2\la)$, we obtain that 
$$K(S) \geq \left( \sigma^2 - \la \right)^2 /(\sigma^2 \la).$$ 

In Example~1 where $S = S_n=\sum_{i=1}^n X_i$
is the sum of $n$ $\iid$ Bernoulli($\la/n$) random
variables, the lower bound (\ref{eq:CR}) coincides 
with the upper bound given in Proposition~3. 
Similarly, in Example~3 with all the $q_i=n/(n+\la)$, 
the upper bound from Proposition~3 holds
with equality. Therefore, any remaining slackness 
in our bounds comes from either Proposition~2 or Pinsker's
inequality.

\medskip

Finally, in Proposition~4 below we establish a formal connection 
between relative entropy and the probability distribution 
$(x+1)P(x+1)/\la$ implicitly used in our definition 
of the scaled Fisher information. It is proved
in the next section.

\medskip

\noindent
{\em Proposition 4.}
Let $X$ be an integer-valued random variable with
distribution $P$ and mean $\la$. If $X$
is the sum of independent Bernoulli 
random variables, then
\be
D(P \| \Pol) = \int_{0}^{\infty} D(P_t \| \widetilde{P}_t) dt,
\label{eq:dB}
\ee
where $P_t(r)=\Pr(X_t=r)$ is the distribution of 
$X_t=X + \Pot$ where $\Pot$ is an independent
Poisson($t$) random variable, and 
$\widetilde{P}_t(r) = (r+1)\PR(X_t = r+1)/(\la+t)$.
More generally, the same result holds for any random
variable $X$ that has $K(X)<\infty$ and satisfies 
the logarithmic Sobolev inequality of Proposition~2.

\medskip

This result is reminiscent of the well-known de Bruijn 
identity, which states that the (differential) relative entropy 
between a random variable $X$ and a Gaussian with the same
variance can be written as a weighted
integral of (continuous) Fisher informations of convex
combinations of $X$ and an independent $N(0,t)$ random
variable; see \cite{cover:book}\cite{barron:clt}.
In a similar vain,
if we formally expand the logarithm in the
integrand in (\ref{eq:dB}) as a Taylor
series, then the first term in the expansion
(the quadratic term) turns out to be equal to
$K(X_t)/2(\lambda+t)$. Therefore, 
$$D(P \| \Pol) \approx \int_{0}^{\infty} 
\frac{K(X+\Pot)}{2(\lambda+t)} dt,$$ 
giving an alternative formula 
to Proposition~2, also relating 
scaled Fisher information and relative
entropy.


\subsection{Proofs}

Although in several places below we formally
divide by a quantity which may be zero,
this is taken care of by the usual
conventions, $0\log (0/a)=0$,
$0\log (0/0)=0$, and
$0\log (a/0)=\infty$,
for any $a>0$.

\medskip

\noindent
{\em Proof of Proposition~2. }
Let $\mbox{Po}_\la(k)$ denote the $\Pol$ probabilities.
In the case when $P$ has full support, the result follows
immediately from Corollary~4 of \cite{bobkov-ledoux:98},
upon considering the function 
$f(k)=P(k)/\mbox{Po}_\la(k),$ $k\geq 0$.

In the case of finite support, for $\epsilon>0$ let
$X^\epsilon$ have the mixture distribution
$$P^\epsilon = (1-\epsilon)\mbox{Po}_\la +\epsilon P.$$
Then $E(X^\epsilon)=\la$ and $P^\epsilon$ has full
support, so by the previous part,
\be
D(P^\epsilon\|\Pol)\leq K(X^\epsilon).
\label{eq:full}
\ee
But since $P(k)=0$ for $k\geq N+1$, then
$P^{\epsilon}(k)/P(k)=\epsilon$ for those $k$,
and letting $\epsilon\downarrow 0$ in the
left hand side of (\ref{eq:full}) we get
$$D(P^\epsilon\|\Pol)=\sum_{k=0}^NP^{\epsilon}(k)\log\Big[
\frac{P^{\epsilon}(k)}{\mbox{Po}_\la(k)}\Big]
+\PR\{\Pol>N\}\epsilon
	\log\epsilon\to D(P\|\Pol).$$
Moreover, 
$$\frac{(k+1)P^\epsilon(k+1)}{\la P^{\epsilon}(k)}=1,
\;\;\;\;k\geq N+1,$$
so
$$K(X^\epsilon)=\sum_{k=0}^NP^{\epsilon}(k)
\Big[\frac{(k+1)P^\epsilon(k+1)}{\la P^{\epsilon}(k)}-1
\Big]^2\to K(X),$$
as $\epsilon\downarrow 0$, and this completes the proof.
\qed

\medskip

Next we prove the bound in (\ref{eq:tvbound}) using
a classical Poincar\'{e} inequality for the Poisson
distribution. We actually establish the following 
(apparently stronger) bound
for the Hellinger distance $\| P- \Pol \|_H$ 
between $P$ and $\Pol$:
$$ \|P-\Pol\|^2_\stv \leq \| P-\Pol\|^2_H \leq 2 K(X).$$

\noindent
{\em Proof of (\ref{eq:tvbound}). }
For any function $f: \IN_+ \rightarrow \RL$, define
$\Delta f (x) = f(x+1) - f(x)$.
It is well-known that,
writing $\mbox{Po}_\la(x)$ for the Poisson($\lambda$) probabilities,
then for all functions $g$ in $L^2(q)$,
\be
\sum_x \mbox{Po}_\la(x) (g(x) - \mu)^2 
\leq \lambda \sum_x \mbox{Po}_\la(x) (\Delta g(x))^2,
\label{eq:poincare}
\ee
where $\mu = \sum_x g(x) \mbox{Po}_\la(x)$
is the mean of $g$ under $\Pol$;
see for example Klaassen \cite{klaassen:85}.

Using the simple fact that
$$(\sqrt{u} -1)^2 \leq (\sqrt{u} -1)^2
(\sqrt{u} +1)^2 = (u-1)^2, \;\;\;\;\mbox{for all $u\geq 0$},$$
we get that
$$K(X)
=
\lambda \sum_x P(x) \left( \frac{ P(x+1) \mbox{Po}_\la(x)}
{\mbox{Po}_\la(x+1) P(x)} - 1
        \right)^2
\geq \lambda \sum_x P(x) \left( \sqrt{\frac{ P(x+1) \mbox{Po}_\la(x)}
{\mbox{Po}_\la(x+1) P(x)}}
        - 1 \right)^2,
$$
and applying~(\ref{eq:poincare})
to the function $g(x) =
\sqrt{P(x)/\mbox{Po}_\la(x)}$
we obtain,
\begin{eqnarray*}
K(X)
& \geq & \lambda \sum_x \mbox{Po}_\la(x) 
\left( \sqrt{\frac{ P(x+1)}{\mbox{Po}_\la(x+1)}}
- \sqrt{\frac{ P(x)}{\mbox{Po}_\la(x)}} \right)^2  \\
& \geq
& \sum_x \mbox{Po}_\la(x) 
\left( \sqrt{\frac{ P(x)}{\mbox{Po}_\la(x)}} - \mu \right)^2
\; =\; 1 - \mu^2,
\end{eqnarray*}
where $\mu = \sum_x \sqrt{P(x)\mbox{Po}_\la(x)}$.
Therefore, the Hellinger distance $\|P-\Pol\|_H$ satisfies
$$\| P-\Pol \|^2_H = (2 - 2\mu) \leq 2(1-\mu^2) \leq 2 K(X),$$
and since
$ \|P-\Pol\|_\stv \leq \sqrt{ \| P- \Pol \|_H }$
(see, e.g., \cite[p.~360]{saloff-coste:97})
the result follows.
\qed

\medskip

For the proof of Proposition~3,
as in the case of normal convergence in Fisher information, 
we exploit the theory of $L^2$ spaces and the fact that 
scaled score functions of sums are 
conditional expectations (projections) 
of the original scaled score functions.

\medskip

\noindent
{\em Lemma. Convolution}:
If $X$ and $Y$ are nonnegative integer-valued
random variables with probability 
distributions $P$ and $Q$ and means $p$ and $q$, respectively,
then,
$$ \ro_{X+Y}(z) = 
E [\alpha_X \ro_X(X) + \alpha_Y \ro_Y(Y) \mid X+Y = z] ,$$
where $\alpha_X = p/(p+q)$, $\alpha_Y = q/(p+q)$.

\smallskip

\noindent
{\em Proof. }
Writing $F(z+1) = \sum_x P(x) Q(z-x+1)$ for 
the distribution of $X+Y$, we have,
\begin{eqnarray*}
\ro_{X+Y}(z) 
& = & \sum_x \frac{ (z+1) P(x) Q(z-x+1)}{(p+q) F(z)} - 1 \\
& = & \sum_x \Big[ \frac{ x P(x) Q(z-x+1)}{(p+q) F(z)} 
	+  \frac{ (z-x+1) P(x) Q(z-x+1)}{(p+q) F(z)}\Big] - 1 \\
& = & \alpha_X \left[ \sum_x \frac{ x P(x)}{p P(x-1)} 
	\frac{P(x-1) Q(z-x+1)}{ F(z)} - 1 \right] \\
& & + \alpha_Y \left[ \sum_x \frac{ (z-x+1) Q(z-x+1)}{q Q(z-x)} 
	\frac{P(x) Q(z-x)}{F(z)} - 1 \right] \\
&\eqa& \sum_x \frac{P(x) Q(z-x)}{F(z)} 
	\Big[ \alpha_X \ro_X(x) + \alpha_Y \ro_Y(z-x) \Big],
\end{eqnarray*}
as required, where~$(a)$ follows 
by moving $x$ to $(x+1)$ in the first sum.
\qed

\medskip

\noindent
{\em Proof of Proposition~3. }
It suffices to prove the case $n=2$. 
By the Lemma,
\begin{eqnarray*}
0 & \leq & E \left[ \frac{p_1}{\la} \ro_{X_1}(X_1) 
	+\frac{p_2}{\la} \ro_{X_2}(X_2) - \ro_{X_1+X_2}(X_1+X_2)
\right]^2 \\
& = & E \left[\frac{p_1}{\la} \ro_{X_1}(X_1)
        +\frac{p_2}{\la} \ro_{X_2}(X_2) \right]^2 - 
	E [\ro_{X_1+X_2}(X_1+X_2)]^2 ,
\end{eqnarray*}
therefore, noting that $E[\ro_{X}(X)]=0$
for any random variable $X$,
\begin{eqnarray*}
K(X_1+X_2) 
& = & 
	(p_1+p_2) E [\ro_{X_1+X_2}(X_1+X_2)]^2 \\ 
& \leq &
	\la E \left[ \frac{p_1}{\la} \ro_{X_1}(X_1) +\frac{p_2}{\la} 
	\ro_{X_2}(X_2) \right]^2 \\
& = & 
	\frac{p_1}{\la} \Big(p_1 E [\ro_{X_1}(X_1)]^2\Big) + 
	\frac{p_2}{\la} \Big(p_2 E [\ro_{X_2}(X_2)]^2\Big) \\
& = & \frac{p_1}{\la} K(X_1) + \frac{p_2}{\la} K(X_2),
\end{eqnarray*}
as claimed.
\qed

\medskip

\noindent
{\em Proof of Proposition~4. }
Assume for the moment that the relative entropy 
between $P_t$ and Po($\la+t$) tends to zero as $t\to\infty$ 
(this will be established below).
Then we can write 
$D(P\|\Pol)$ as the integral
\begin{eqnarray*}
D(P \| \Pol) & = & - \int_{0}^{\infty}
\frac{\partial}{\partial t} D(P_t \| \mbox{Po}(\lambda+t)) dt \\
& = & - \int_{0}^{\infty} \frac{\partial}{\partial t} 
	\Big(
	(\lambda +t ) - E[X_t \log (\lambda+t)]
	+ E[\log (X_t!)] - H(X_t) 
	\Big) dt \\
& = & 
	\int_{0}^{\infty} \left( \log(\lambda +t ) 
	- \frac{\partial}{\partial t}
	E[\log (X_t!)] +  \frac{\partial}{\partial t} H(X_t) \right) dt.
\end{eqnarray*}
Since the probabilities
$P_{t}$
satisfy a differential-difference equation,
$ \frac{ \partial P_{t}}{\partial t}(x) = P_{t}(x-1) - P_{t}(x),$
we have,
$$\frac{\partial}{\partial t} E [\log (X_t)!] =
\sum_r \frac{\partial P_t}{\partial t}(r) \log r!
= \sum_r (P_t(r-1) - P_t(r))\log r! = E \log(X_t +1),$$
and similarly,
\begin{eqnarray*}
\frac{\partial}{\partial t} H(X_t) & = & - \sum_r
\frac{\partial P_t}{\partial t}(r) \log P_t(r)
=   \sum_r P_t(r) \log \left( \frac{P_t(r)}{P_t(r+1)} \right).
\end{eqnarray*}
Substituting these two expressions 
in the expansion of $D(P\|\Pol)$ 
the result follows.

Finally it remains to establish our initial assumption.
If $X$ is the sum of independent Bernoulli random variables
then it has finite support and Proposition~2 holds; moreover,
$K(X)$ is easily seen to be finite by Proposition~3. More 
generally, using Propositions~2 and~3 we have
\ben
D(P_t\|\mbox{Po}(\la+t))
\leq
	K(X+\Pot)
\leq
	\frac{\la}{\la+t} K(X)\to 0,
\een
as $t\to\infty$, as required.

\qed

\section*{Acknowledgments}
Interesting conversations with
Persi Diaconis, Stu Geman and Matt Harrison
are greatfully acknowledged. We also
thank the Associate Editor and the two
referees for their useful comments on an
earlier version of the manuscript.


\end{document}